\newif \ifdraft
\newif\ifarxiv
\arxivtrue

\documentclass[reqno,  12pt]{article}

\ifdraft
\usepackage[12pt]{extsizes}

\usepackage[inline]{showlabels}
\usepackage[pagewise]{lineno}
\linenumbers
\fi

\usepackage[a4paper,  margin = .75in]{geometry}

\usepackage{amsmath,  amsthm,  amssymb,  mathtools}
\numberwithin{equation}{section}

\usepackage{thmtools}
\usepackage{enumitem}

\usepackage[dvipsnames]{xcolor}
\usepackage{hyperref}
\hypersetup{
	colorlinks  =  true, 
	linkcolor  =  {BlueViolet}, 
	citecolor = {RoyalBlue}, 
	urlcolor = {BrickRed}
}

\makeatletter
\renewcommand{\itemautorefname}{\@gobble}
\makeatother

\newtheorem{theorem}{Theorem}[section]
\newtheorem{proposition}[theorem]{Proposition}
\newtheorem{lemma}[theorem]{Lemma}
\newtheorem{corollary}[theorem]{Corollary}
\theoremstyle{definition}

\newtheorem*{theorem*}{Theorem}
\newtheorem{theoremX}{Theorem}

\newcommand{\R}{\mathbb{R}}
\newcommand{\C}{\mathbb{C}}
\newcommand{\T}{\mathcal{T}}
\newcommand{\Hit}{\mathcal{H}_n}
\newcommand{\Z}{\mathcal{Z}_n}
\newcommand{\Q}{\mathcal{Q}}
\newcommand{\PSL}{\mathrm{PSL}}
\newcommand{\SL}{\mathrm{SL}}
\newcommand{\Crit}{\operatorname{Crit}}
\newcommand{\Hess}{\operatorname{Hess}}
\newcommand{\Area}{\operatorname{Area}}
\newcommand{\ind}{\operatorname{ind}}

\title{Persistent critical points and Morse energies\\
for Hitchin representations}

\author{
	Prahllad Deb\thanks{
		Department of Mathematics,  Indraprastha Institute of Information Technology Delhi,  New Delhi,  India.
		\texttt{\href{mailto:prahllad@iiitd.ac.in}{prahllad@iiitd.ac.in}}
	}
	\and
	Anu Dhochak\thanks{
		International Centre for Theoretical Sciences (ICTS-TIFR),  Bengaluru,  India.
		\texttt{\href{mailto:anu.dhochak@icts.res.in}{anu.dhochak@icts.res.in}}
	}
	\and
	Pradip Kumar\thanks{
		Department of Mathematics,  School of Natural Sciences,  Shiv Nadar Institution of Eminence,  Delhi-NCR,  India.
		\texttt{\href{mailto:pradip.kumar@snu.edu.in}{pradip.kumar@snu.edu.in}}
	}
	\and
	Sachchidanand Prasad\thanks{
		Department of Mathematics,  Indraprastha Institute of Information Technology Delhi,  New Delhi,  India.
		\texttt{\href{mailto:sachchidanand@iiitd.ac.in}{sachchidanand@iiitd.ac.in}}
	}
}

\newcommand{\setsubjclass}[1]{\def\thesubjclass{#1}}
\newcommand{\setkeywords}[1]{\def\thekeywords{#1}}
\newcommand{\printclassification}{%
	\begingroup
\renewcommand{\thefootnote}{}%
\begin{NoHyper}%
	\footnotetext{\textbf{2020 Mathematics Subject Classification.} \thesubjclass.}%
	\footnotetext{\textbf{Keywords.} \thekeywords.}%
	\end{NoHyper}%
\endgroup
}

\setsubjclass{Primary: 53C43,  32G15; Secondary: 53C42,  57R70,  22E40}

\setkeywords{Hitchin representations, equivariant harmonic maps, 
Teichm\"uller space, minimal surfaces,  Morse theory,  generic nondegeneracy}

\begin{document}
\date{}
\maketitle
\printclassification

\vspace{-1em}
\begin{abstract}
We show that the energy on Teichm\"uller space is Morse for a dense $G_\delta$ set of full measure in each Hitchin component of $\PSL(n,\ \R)$. Suppose a representation has an unstable critical point, possibly a degenerate one. We show that every neighbourhood of that representation contains a nonempty open set on which three distinct nondegenerate critical points vary real-analytically with the representation. One is a strict local minimum. Another has positive Morse index and can be chosen arbitrarily close to the original unstable point. Applying this result to the unstable minimal surfaces of Sagman and Smillie gives such families in every Hitchin component when the genus satisfies $g \geq 3$ and $n \geq 4$.
\end{abstract}

\section{Introduction and main results}\label{sec:introduction}
Equivariant minimal surfaces associated to a surface-group representation arise as critical points of an energy function on Teichm\"uller space.
For Hitchin representations. In ~\cite{LabourieEnergy} ~\cite{LabourieRankTwo}, Labourie proves the existence and uniqueness in real rank two.  In real rank at least three, 
Sagman--Smillie prove nonuniqueness in the case of genus at least
three~\cite{SS}.  We study the nondegeneracy
of these critical points and the persistence of their multiplicity under
deformation of the representation.

Let $S$ be a closed oriented surface of genus $g  \, \geq \,  2$,  write
$\T = \T(S)$,  and fix a Hitchin component
$$\Hit \subset
  \mathrm{Hom}_{\mathrm{red}}(\pi_1(S),\ \PSL(n,  \R))/\PSL(n,  \R), 
 \qquad n  \, \geq \,  2.$$
Let  $X_n = \SL(n,\ R)/\mathrm{SO}(n)$ with a fixed invariant metric.
For $\rho  \in  \Hit$ and $J  \in  \T$,  let
$u_{\rho,  J}  \colon  \widetilde S \to X_n$ be the unique
$\rho$-equivariant harmonic map, and define
$$E_\rho(J)  =  \frac12  \int_S \|du_{\rho,  J}\|_{g_J}^2 \,  dA_{g_J}, $$
where $g_J$ is any metric in the conformal class $J$.   The energy is independent of the metric in the conformal class.  The first variation of the function $E_\rho$ is zero if and only if the Hopf
differential of $u_{\rho, J}$ vanishes.  Thus its critical points determine equivariant minimal immersions
\cite[Remark~4.3]{SS}. Throughout, by the distinct critical points, we mean distinct marked conformal structures in $\T(S)$.

A critical point is \emph{nondegenerate} if the kernel of its Hessian is trivial, and it is called \emph{unstable} if the Hessian has a negative direction. Its \emph{Morse index} is the maximal dimension of a real subspace on which the Hessian is negative definite. The energy is \emph{Morse} if
all its critical points are nondegenerate.

\medskip 

We prove the following:
\begin{theoremX}\label{intro:morse}
Let $g\, \geq \, 2$ and $n\, \geq \, 2$. The set
$$\mathcal R_n\,  = \, 
 \bigl\{\rho\,  \in \, \Hit\, \big|\, 
 E_\rho\text{ is Morse on }\T\bigr\}$$
is a dense $G_\delta$ subset of $\Hit$ and has full measure in
every coordinate chart. For each $\rho\, \in \, \mathcal R_n$
and each $B \,\in \, \R$,  $\bigl\{J\,  \in \, \T\, \big|\, 
E_\rho(J)\, \leq\, B\bigr\}$ contains only finitely many
critical points.
\end{theoremX}

\begin{theoremX}\label{intro:persistence}
Let $g\, \geq \, 2$ and $n\, \geq \, 2$,  let $\rho_*\,  \in \, \Hit$,  and suppose
$J_*\,  \in \, \Crit(E_{\rho_*})$. Assume that
$\Hess_{J_*}E_{\rho_*}$ has a negative direction.  Then for every open neighbourhood $V\, \subset\, \Hit$
of $\rho_*$ and every open neighbourhood $O\, \subset\, \T$ of $J_*$, 
there exist a nonempty connected open set $W\, \subset\, V$ and
real-analytic maps
$$ J_0, \, J_1, \, J_2 \colon  W\,  \longrightarrow \, \T$$
such that,  for every $\rho\,  \in \, W$,  their values are pairwise
distinct nondegenerate critical points of $E_\rho$. 
$J_0(\rho)$ is a strict local minimum.  Moreover,  $J_1(W)\, \subset\, O$, 
and $J_1(\rho)$ has constant Morse index $k$ on $W$,  where
$1\, \leq\, k\, \leq\, 3g-3$. If the Hessian is negative definite
on an $r$-dimensional real subspace, then the branch $J_1$ can be
chosen with $k\, \geq \, r$.
\end{theoremX}
The open set $W$ need not contain $\rho_*$, and only the unstable branch
is required to stay near $J_*$.  We have a corollary~\ref{cor:open-families}, 
if $g \, \geq \, 3$ and $n \, \geq \, 4$, every Hitchin component contains a nonempty
open set carrying three distinct nondegenerate real-analytic marked
critical branches, one consisting of strict local minima and one of
constant positive Morse index. Here the unstable examples come
from~\cite[Theorem~A]{SS}; the area-to-energy comparison is recalled in
Lemma~\ref{lem:area}.

We also give bounded-energy continuation as
Proposition~\ref{prop:bounded-continuation}. As an application
of~\cite[Definition~4.11 and Lemma~4.10]{SS}, we see that near a Morse representation, 
all critical points below a fixed noncritical energy level continue in finitely many nondegenerate real-analytic branches of constant index.

\smallskip

We use Slegers's theorem on the real-analytic dependence of harmonic maps~\cite{SlegersAnalytic}, together with the joint moduli space of Higgs bundles and the relative Hitchin map constructed by Collier, Toulisse, and Wentworth~\cite{CTW}.

Section~\ref{sec:preliminaries} introduces the notation. In Section~\ref{sec:incidence}, we study the space of critical pairs and use it to prove Theorem~\ref{intro:morse} on generic nondegeneracy. In Section~\ref{sec:stable}, we prove Theorem~\ref{intro:persistence} on the persistence of critical points and derive its consequences.
\section*{AI disclaimer} No LLM was used to prepare this manuscript.  Copilot (built-in tool in LaTeX editor) was used as editorial assistant to improve sentences and grammar. 
\section{Preliminaries}\label{sec:preliminaries}

For $J\,  \in \, \T$,  let $K_J$ denote the holomorphic cotangent bundle of
$(S, J)$. For fixed $J$, Hitchin's parametrisation is the real-analytic
diffeomorphism
\begin{equation}\label{eq:hitchinintro}
 h_J\  \colon \ \bigoplus_{a\,  = \, 2}^n H^0(S, \, K_J^a)\,  \longrightarrow \, \Hit, 
 \qquad(q_2, \, \ldots, \, q_n)\, \longmapsto\, \rho.
\end{equation}
Here $q_a$ is a holomorphic differential of degree $a$,  and $h_J$ is
the Hitchin section followed by nonabelian Hodge theory and
holonomy~\cite{Hitchin, Corlette, Donaldson}.  For odd $n$, there is a unique Hitchin component
of $\PSL(n,\ R)$,  whereas for even $n$ there are
two. In the latter case, we choose the Hitchin
section on a fixed component $\Hit$, as the other component is obtained by conjugation with a matrix in $\mathrm{GL}(n,\ R)$ of negative
determinant~\cite{Hitchin}.

Write 
\begin{equation}\label{eq:quadratic_bundle}\Q\,  \longrightarrow  \, \T\end{equation} for the holomorphic vector bundle with fibre
$\Q_J\,  = \, H^0(S, \, K_J^2)$. Its complex rank is $m\,  = \, 3g-3$. Since
$u_{\rho, J}$ is harmonic, its Hopf differential is holomorphic. It
therefore defines a section of the pullback of $\Q$ to $\Hit \times \T$:
$$ \Phi(\rho, \, J)
 \, : = \, \langle {(u_{\rho, J})}_z, \, {(u_{\rho, J})}_z\rangle\, dz^2
 \,  \in \, \Q_J.$$
 We use the normalisation in which, for a Beltrami differential
$\mu$ representing a tangent vector to $\T$, 
\begin{equation}\label{eq:variation}
 dE_\rho\big|_J(\mu)\,  = \, -4\operatorname{Re}
  \int_S\Phi(\rho, \, J)\, \mu.
\end{equation}

Its real part is nondegenerate and gives a real-linear
isomorphism from the underlying real vector space of
$H^0(S, \, K_J^2)$ to $T_J^*\T$. We have
\begin{equation}\label{eq:zero}
                     dE_\rho\big|_J\,  = \, 0\quad\Longleftrightarrow\quad
                     \Phi(\rho, \, J)\,  = \, 0.
\end{equation}
Thus,  $\Phi\,  = \, 0$ is a global equation for the critical pairs. Set
\begin{equation}\label{eq:incidence}
 \Z\,  = \, \{(\rho, \ J) \in \Hit \times \T:dE_\rho(J) = 0\}, 
 \qquad p(\rho, J)\,  = \, \rho.
\end{equation}

We also use the holomorphic vector bundle
\begin{equation}\label{eq:minbundle}
 \mathcal B_n\,  \longrightarrow \, \T, \qquad
 (\mathcal B_n)_J\,  = \, \bigoplus_{a\,  = \, 3}^n H^0(S, \, K_J^a),  
\end{equation}
with the convention that $\mathcal B_2$ is
the zero vector bundle over $\T$.
The holomorphic structure on $\mathcal B_n$
is induced by the direct images of powers
of the relative canonical bundle over the
universal Teichm\"uller curve. Indeed,  by
Riemann--Roch, 
\begin{equation}\label{eq:dimH_0(S, K_J^a)}
\dim_{\mathbb C}H^0(S, K_J^a)
 = (2a-1)(g-1),  \qquad a \, \geq \,  3, 
\end{equation}
independently of $J$. Thus Grauert's theorem
shows that these spaces form holomorphic
vector bundles over $\T$.
Its fibre is zero when $n\,  = \, 2$. The associated map is
\begin{equation}\label{eq:labourieintro}
 \mathcal B_n\,  \longrightarrow \, \Hit, \qquad
 (J, \, q_3, \, \ldots, \, q_n)\, \longmapsto\, 
 h_J(0, \, q_3, \, \ldots, \, q_n).
\end{equation}
\section{The space of critical pairs}\label{sec:incidence}
\subsection{Analytic dependence and the energy}
\begin{proposition}\label{prop:analytic-dependence}
\setlength{\emergencystretch}{1em}
The function $(\rho, J)\longmapsto E_\rho(J)$
is jointly continuous  and the Hopf section
$(\rho, J)\longmapsto\Phi(\rho, J)$
is jointly real analytic on $\Hit \times \T$.
For each fixed $\rho$,  the function
$J\longmapsto E_\rho(J)$ is real analytic.

\end{proposition}

\begin{proof}
Let $\Gamma\,  = \, \pi_1(S)$,  and let
$\operatorname{Hom}_{\mathrm{Hit}}(\Gamma, \, \PSL(n, \, \R))$ denote actual
Hitchin representations, before quotienting by conjugation. From~\cite{SlegersAnalytic} it follows that, for every $k \in \mathbb N$ and $0 < \alpha < 1$, the map
$$
 (J, \, \widetilde\rho)\, \longmapsto\, u_{\widetilde\rho, J}
$$
is continuous with values in
$C^{k, \alpha}(\widetilde S, \, X_n)$. If
$\Omega\, \subset\, \widetilde S$ is a bounded domain containing a
fundamental domain, the restricted map
$$
 (J, \, \widetilde\rho)\, \longmapsto\, 
 u_{\widetilde\rho, J}|_\Omega
$$
is real analytic with values in
$C_b^{k, \alpha}(\Omega, \, X_n)$
as stated in~\cite[Proposition 2.12]{SlegersAnalytic}.

Use the real-analytic family of representative metrics $g_J$ as in \cite{SlegersAnalytic}. For $k\, \geq \, 1$,  forming the
energy density and integrating it over a fundamental domain show that
the energy is continuous in both variables and real analytic in $J$
for each fixed homomorphism. It is invariant under conjugation, so it
descends continuously to $\Hit \times \T$ and remains real analytic in
$J$ for each fixed conjugacy class.

For joint real analyticity of the Hopf section, set
$G_{\C} = \PSL(n, \C)$. By~\cite[Lemma~10.1 and
Theorem~1.5]{LabourieAnosov},  every Hitchin representation
is irreducible over $\R$,  and each nontrivial element has
a lift to $\SL(n,\ R)$ with distinct real eigenvalues.
These properties imply irreducibility over $\C$.
Indeed, let $W\subset\C^n$ be a complex subspace invariant
under the full inverse image of $\rho(\Gamma)$ in
$\SL(n,\ R)$. Choose an element of this inverse image
with distinct real eigenvalues. In a real eigenbasis
for this element,  $W$ is the complex span of a subset
of the basis vectors. Their real span $V$ satisfies
$$ W = V+iV,  \qquad V = W \cap \R^n. $$
Since the lifted group consists of real matrices and
preserves $W$, it also preserves $V$. Real irreducibility
therefore gives $V = \{0\}$ or $V = \R^n$,  and hence
$W = \{0\}$ or $W = \C^n$. Therefore, the inclusion $\PSL(n,\R)$ in $\PSL(n,\C)$ induces a real-analytic map of character spaces 
\begin{equation}\label{eq:icrho}
\iota_{\C} \colon \Hit \longrightarrow 
\mathcal X^{\mathrm{irr}}(G_{\C}).
\end{equation}

Let $\mathbf M(G_{\C})$ be the joint moduli space of stable
$G_{\C}$-Higgs bundles constructed in~\cite{CTW},  with
holomorphic projection $\pi$ and nonabelian Hodge map
$\mathcal H$. By~\cite[Section~4.9 and
Theorem~4.23]{CTW},  the map
\begin{equation}\label{eq:pi_H}
F = (\pi, \mathcal H) \colon 
\mathbf M(G_{\C})
 \longrightarrow 
\T \times \mathcal X^{\mathrm{irr}}(G_{\C})
\end{equation}
is a real-analytic homeomorphism, and for each $J \in \T$
its restriction $\mathcal H_J \colon \pi^{-1}(J)
 \longrightarrow \mathcal X^{\mathrm{irr}}(G_{\C})$
is a real-analytic diffeomorphism of orbifolds.
In compatible charts on the domain and
codomain of $F$,  its differential $dF$ is invertible. Using the real-analytic inverse
function, we prove that $F^{-1}$ is real analytic.

Let $\varpi$ denote the relative Hitchin map \cite[Theorem~7.3]{CTW},  and let
$$\varpi^{(2)} \colon \mathbf M(G_{\C}) \longrightarrow \Q$$
be its holomorphic quadratic differential.
By~\cite{CTW}, this is a fixed nonzero scalar multiple of the Hopf differential of the corresponding harmonic
map into $G_{\C}/\mathrm{PSU}(n)$.

To compare this differential with $\Phi$, consider
the natural totally geodesic embedding
$$\jmath \colon  X_n \longrightarrow  G_{\C}/\mathrm{PSU}(n).$$
The composition $\jmath\circ u_{\rho, J}$ is harmonic
and equivariant for the complexified representation.    Since the harmonic map is unique, it is the harmonic map corresponding
to $F^{-1}(J, \iota_{\C}(\rho))$. Restricting the metric of the complex symmetric space to the $X_n$, turns out to the invariant metric of $X_n$, and hence  there is a constant
$c \in \C \setminus \{0\}$ such that
$$ \Phi(\rho, J) =  c \, \varpi^{(2)}\bigl(F^{-1}(J, \iota_{\C}(\rho))\bigr).
$$
The right-hand side is a composition of real-analytic
maps. Therefore $\Phi$ is jointly real analytic.
\end{proof}

Recall that $\Q_{\R}$ denotes the underlying real
vector bundle of $\Q$ (see~\eqref{eq:quadratic_bundle}).
For each $J \in \T$,   since the pairing in~\eqref{eq:variation} is non-degenerate, hence we get a linear isomorphism of vector spaces
$$ A_J \colon \Q_J \longrightarrow  T_J^*\T$$ defined by$A_J(q)(v) = -4\operatorname{Re} \int_S q\, \mu, $
where $\mu$ is a tangent vector $T_J\T$. This gives the isomorphism of real vector bundles
$$ A \colon \Q_{\R} \longrightarrow  T^*\T. $$

By~\eqref{eq:variation},  we have
$$ dE_\rho\big|_J = A_J(\Phi(\rho, J)). $$
At a critical pair $(\rho, J) \in \Z$,  we have
$\Phi(\rho, J) = 0$. Differentiating this identity with respect to $J$ and applying the product
rule gives 
\begin{equation}\label{eq:hessian}
\Hess_J E_\rho
 = A_J\circ D_J\Phi(\rho, J), 
\end{equation}
because the term involving the derivative
of $A_J$ vanishes at $\Phi(\rho, J) = 0$.

Here the Hessian is viewed as a linear map
$T_J\T \longrightarrow  T_J^*\T$.
At a zero of $\Phi$,  the derivative $D_J\Phi$
is independent of the chosen local trivialisation
of $\Q$. Consequently,  the Hessian depends
continuously on $(\rho, J) \in \Z$.
\subsection{Hitchin coordinates and transversality}\label{subsec:transversality}
The Hopf section is transverse to the zero section.  Recall the Hitchin
parametrisation~\eqref{eq:hitchinintro},  where
$q_2$ is the quadratic Hitchin coordinate.
By~\cite[Section~7.1.3 and Proposition~7.8]{CTW},  further,  the degree-two component of the Hitchin map
is a fixed nonzero scalar multiple of the
Hopf differential. Since every degree-two
Ad-invariant polynomial on $\mathfrak{sl}(n, \C)$
is a scalar multiple of $X\mapsto\operatorname{tr}(X^2)$, 
there exists a constant $c_n \in \C \setminus \{0\}$, 
independent of $J$ and of the Hitchin-base point, 
such that
\begin{equation}\label{eq:hopf-q2}
\Phi\bigl(h_J(q_2, \ldots, q_n), J\bigr)
 = c_nq_2.
\end{equation}
This constant depends only on the fixed
normalisations. 

\begin{lemma}\label{lem:transverse}
For each fixed $J \in \T$ and every $\rho \in \Hit$, 
the derivative
$D_\rho\Phi(\rho, J):
T_\rho \Hit \longrightarrow \Q_J$ 
is surjective. Therefore, the Hopf section
$\Phi$ is transverse to the zero section.
\end{lemma}

\begin{proof}
Fix $(\rho, J) \in \Hit \times \T$ and write
$\mathbf q = h_J^{-1}(\rho)$. Consider the 
inclusion 
$$ \iota_2 \colon \Q_J \longrightarrow 
\bigoplus_{a = 2}^n H^0(S, K_J^a), 
\qquad q\longmapsto(q, 0, \ldots, 0). $$
Differentiating~\eqref{eq:hopf-q2} in the
$q_2$-direction,  while keeping $J$ and
$q_3, \ldots, q_n$ fixed,  gives
$$
D_\rho\Phi(\rho, J)\circ(dh_J)_{\mathbf q}
\circ\iota_2
 = c_n\operatorname{id}_{\Q_J}.
$$
Since $c_n\neq0$, the right-hand side is
an isomorphism. Therefore $D_\rho\Phi(\rho, J)$
is surjective.

At every zero of $\Phi$, surjectivity
of the derivative in the representation
direction implies surjectivity of the
full derivative onto $\Q_J$. Hence
$\Phi$ is transverse to the zero section.
\end{proof}
\subsection{The critical-pair space and its projection}\label{subsec:incidence-manifold}

Now we return to the critical-pair space and its projection as in 
in~\eqref{eq:incidence}.  We have $\Z\,  = \, \{\Phi\,  = \, 0\}$ ~\eqref{eq:zero}, and   Proposition~\ref{prop:analytic-dependence}, 
Lemma~\ref{lem:transverse}, together with the  implicit function theorem
show that $\Z$ is a closed real-analytic submanifold of
$\Hit   \times  \T$.  Since $\Phi$ is transverse to the zero section, 
$\Z$ is a real-analytic submanifold of
$\Hit \times \T$ of real codimension $6g-6$.
As $\dim_\R\T = 6g-6$,   and from \eqref{eq:dimH_0(S, K_J^a)},  we have
$$
\dim_\R\Z = \dim_\R\Hit \  =  \ 2(n^2-1)(g-1).
$$

Now we relate the differential of 
$p \colon \Z \longrightarrow \Hit$ to the Hessian
of the energy function $E_\rho$ at its critical points.
\begin{lemma}\label{lem:regular}
At $(\rho, \, J)\,  \in \, \Z$ there is a natural identification
$$\ker\bigl(d p_{(\rho, \, J)}\bigr) \, \cong\, \ker\bigl(\Hess_J E_\rho\bigr).
$$
The isomorphism implies that $(\rho, \, J)$ is a regular point of $p$ if and only if $J$ is a
nondegenerate critical point of $E_\rho$.
\end{lemma}

\begin{proof}
Since $\Z = \Phi^{-1}(0)$ and $\Phi$ is
transverse to the zero section, the tangent
space of $\Z$ is the kernel of the
linearisation of $\Phi$ at each zero.
Writing a tangent vector as
$(\dot\rho, \dot J) \in 
T_\rho\Hit\oplus T_J\T$,  we obtain
$$
T_{(\rho, J)}\Z
 = 
\left\{
(\dot\rho, \dot J) \in 
T_\rho\Hit\oplus T_J\T
\, \middle|\, 
D_\rho\Phi(\rho, J)(\dot\rho)
+
D_J\Phi(\rho, J)(\dot J) = 0
\right\}.
$$
Since $p(\rho, J) = \rho$,  its differential
satisfies
$dp_{(\rho, J)}(\dot\rho, \dot J) = \dot\rho.$ Therefore, using the tangent-space
description above, we obtain
$$
\ker dp_{(\rho, J)}
 = 
\{(0, \dot J):D_J\Phi(\dot J) = 0\}.
$$
 Since $A_J$ in ~\eqref{eq:hessian} is an isomorphism,  the map $(0, \, \dot J)\, \longmapsto\, \dot J$ gives the stated identification.

 Since $\dim_\R\Z = \dim_\R\Hit$, the
differential $dp_{(\rho, J)}$ is a linear
map between vector spaces of equal dimension.
Hence it is surjective if and only if its
kernel is trivial. By the preceding
identification, this is equivalent to
$\ker\Hess_JE_\rho = \{0\}$.
Therefore,  $(\rho, J)$ is a regular point
of $p$ if and only if $J$ is a
nondegenerate critical point of $E_\rho$.
\end{proof}
The preceding construction gives the local analytic structure of
$\Z$. Its global topology comes from the varying Hitchin
parametrisation. Recall the holomorphic vector bundle
$\mathcal B_n \longrightarrow \T$
defined in~\eqref{eq:minbundle}.

\begin{lemma}\label{lem:global-parametrisation}
The map
\begin{equation}\label{eq:bundlehomeo}
 \Psi \colon \mathcal B_n\,  \longrightarrow \, \Z, \qquad
 (J, \, q_3, \, \ldots, \, q_n)\, \longmapsto\, 
 \bigl(h_J(0, \, q_3, \, \ldots, \, q_n), \, J\bigr)
\end{equation}
is a homeomorphism and therefore in particular $\Z$ is connected.
\end{lemma}

\begin{proof}
For each fixed $J$,  the diffeomorphism $h_J$ and
identity~\eqref{eq:hopf-q2} show that $\Psi$ is a bijection onto
$\Z$. We verify that its inverse is continuous. For
$(\rho, \, J)\,  \in \, \Z$,   $\iota_{\C}(\rho) \in  \mathcal X^{\mathrm{irr}}(\mathrm{PSL}(n,  \C))$ (see \eqref{eq:icrho}).  Let $(\pi,  \  \mathcal{H})$ as in \eqref{eq:pi_H}.   Set
$$
 x\,  = \, {(\pi, \, \mathcal H)}^{-1}\bigl(J, \, \iota_{\C}(\rho)\bigr)
$$
in the joint Higgs-bundle moduli space used in the proof of
Proposition~\ref{prop:analytic-dependence},  and define
$$
 \Theta(\rho, \, J)\,  = \, 
 \bigl(J, \, \varpi^{(3)}(x), \, \ldots, \, \varpi^{(n)}(x)\bigr).
$$
The homeomorphism stated in~\cite[Section 4.9]{CTW} and the holomorphic relative Hitchin map
$\varpi$ of~\cite[Theorem 7.3]{CTW} show that $\Theta$ is continuous and the equation $\Phi(\rho, \, J)\,  = \, 0$
shows that the omitted quadratic coordinate is zero. Hence $\Theta$
takes values in $\mathcal B_n$ (see \eqref{eq:minbundle}). The fixed-surface Hitchin
parametrisation shows that $\Theta$ is a bijection with inverse $\Psi$.

Finally,  $\T$ and every
vector-space fibre of $\mathcal B_n$ are connected,  hence so are
$\mathcal B_n$ and $\Z$.
\end{proof}

Under~\eqref{eq:bundlehomeo},  the projection $p$ is precisely the
Labourie map introduced in~\eqref{eq:labourieintro}. Each $E_\rho$ is
nonnegative and proper by~\cite[Theorem 1.0.3]{LabourieEnergy}, and
therefore attains a global minimum. Equation~\eqref{eq:zero} then shows that
$p$ is surjective.

\subsection{Regular values and Morse energies}\label{subsec:morse-energies}
We prove the following:
\begin{theorem}[Theorem~\ref{intro:morse}]\label{thm:morse-energies}
Let $g\, \geq \, 2$ and $n\, \geq \, 2$. The set
$$
 \mathcal R_n\,  = \, 
 \bigl\{\rho\,  \in \, \Hit\, \big|\, 
 E_\rho\text{ is Morse on }\T\bigr\}
$$
is a dense $G_\delta$ subset of $\Hit$ and has full measure in
every smooth coordinate chart. For each $\rho\,  \in \, \mathcal R_n$
and each $B\,  \in \, \R$,  the sublevel
$\bigl\{J\,  \in \, \T\, \big|\, 
E_\rho(J)\, \leq\, B\bigr\}$ contains only finitely many
critical points.
\end{theorem}
\begin{proof}
By Lemma~\ref{lem:regular},  a representation $\rho$ is a regular value
of $p$ if and only if every critical point of $E_\rho$ is
nondegenerate. Thus by  the definition of $\mathcal R_n$, the regular-value set of $p$ is
$\mathcal R_n$. Let
$C :=\{(\rho, J): \mathrm{ker} dp_{(\rho, J)} \neq \{0\}\}$,  and choose a
countable exhaustion of $\Z$ by compact sets $K_a$.  By Sard's theorem~\cite{Sard},  $p(C)$ has Lebesgue measure zero in every 
coordinate chart.  Further, each $p(C  \cap  K_a)$ is compact and,  being null,  has empty interior. Hence it is closed and nowhere dense. Since
$$ p(C)\,  = \, \bigcup_{a\,  = \, 1}^ \infty p(C \cap  K_a), 
$$
and $\Hit$ is a Baire manifold,  its complement is a dense
$G_\delta$ set.   $p(C)$ has measure zero in every
smooth coordinate chart. This prove first two. 

For a regular value $\rho$,  every critical point is isolated by the
implicit function theorem. The critical set is closed,  and
$\bigl\{E_\rho\, \leq\, B\bigr\}$ is compact for every real $B$. Infinitely many
critical points in the compact set  has  an accumulation point, and which is  
itself a critical point.  This contradicts the fact that every point is isolated. Thus there are only finitely
many. 
\end{proof}

\begin{lemma}\label{lem:analyticdensity}
Let $F \colon  A\,  \longrightarrow \, B$ be a surjective real-analytic map
between real-analytic manifolds of the same dimension.
If $A$ is connected, then the regular points of $F$ form an open dense
subset of $A$.
\end{lemma}

\begin{proof}
Let $d\,  = \, \dim A\,  = \, \dim B$. By Sard's theorem and surjectivity, 
$F$ has a regular point. In local real-analytic coordinates, at the regular point the Jacobian is non-singular, so the set of
regular points is open.

To prove density, we will prove by contradiction.  Suppose a nonempty open subset $U\, \subseteq\, A$
contains only critical points. Therefore, in local coordinates the Jacobian on $U$ is singular. The real-analytic identity theorem and the
connectedness of $A$ imply that it vanishes everywhere.
\end{proof}

We apply the above lemma to the map  \(p\). It shows that  the regular points of \(p\),  equivalently,  the pairs \((\rho, J) \in \Z\) for which \(J\) is a nondegenerate critical point of \(E_\rho\),  are dense in \(\Z\).

\begin{proposition}\label{prop:local-regularisation}
Let $z_*\,  = \, (\rho_*, \, J_*)\,  \in \, \Z$,  and suppose that
$\Hess_{J_*}E_{\rho_*}$ is negative definite on an $r$-dimensional
real subspace,  where $r\, \geq \, 1$. For every neighbourhood
$V\, \subset\, \Hit$ of $\rho_*$ and every neighbourhood
$O\, \subset\, \T$ of $J_*$,  there are a regular point
$\widehat z\,  = \, (\widehat\rho, \, \widehat J)\,  \in \, \Z \cap (V  \times  O)$, 
a connected open set $V'\, \subset\, V$ containing $\widehat\rho$,  and
a real-analytic map $J_1 \colon  V'\,  \longrightarrow \, O$ such that
$(\rho, \, J_1(\rho))\,  \in \, \Z$ for all $\rho\,  \in \, V'$. These critical
points are nondegenerate and have a constant Morse index
$k\, \geq \, r$. Moreover,  $V'$ contains a representation
$\rho^\dagger\,  \in \, \mathcal R_n$.
\end{proposition}
\begin{proof}
Let $L_*\ \subset \ T_{J_*}\T$ be the given $r$-dimensional
subspace on which $\Hess_{J_*}E_{\rho_*}$ is negative definite.
Choose a basis $v_1, \ldots, v_r$ of $L_*$. Using a local
trivialisation of $T\T$ near $J_*$,  extend these vectors to
linearly independent local vector fields $v_1(J), \ldots, v_r(J)$.

By Proposition~\ref{prop:analytic-dependence} and~\eqref{eq:hessian}, 
the symmetric matrix
$$
 \left(
 \Hess_J E_\rho\bigl(v_i(J), v_j(J)\bigr)
 \right)_{i, j = 1}^r
$$
depends continuously on the critical pair $(\rho, J)$.
Since this matrix is negative definite at $z_* = (\rho_*, J_*)$, 
it remains negative definite on a sufficiently small relatively
open neighbourhood
$$
 z_* \in  U\ \subset \ \Z \cap (V  \times  O).
$$
Thus,  at every pair in $U$,  the Hessian is negative definite
on an $r$-dimensional subspace.

By Lemma~\ref{lem:analyticdensity},  $U$ contains a regular point
$\widehat z = (\widehat\rho, \widehat J)$ of $p$.
By  inverse function theorem, we have open
neighbourhoods $U'\ \subset \ U$ of $\widehat z$ and
$V'\ \subset \ V$ of $\widehat\rho$ such that
$$
 p|_{U'} \colon  U' \longrightarrow  V'
$$
is a real-analytic diffeomorphism. Since $p$ is the projection
onto the first factor, its inverse has the form
$$
 \rho\longmapsto \bigl(\rho, J_1(\rho)\bigr)
$$
for a real-analytic map $J_1 \colon  V' \longrightarrow  O$.
Every pair on this graph belongs to $\Z$,  so $J_1(\rho)$
is a critical point of $E_\rho$.

The differential of $p$ is invertible throughout $U'$.
Hence Lemma~\ref{lem:regular} implies that all these critical
points are nondegenerate. Shrinking $V'$ to a connected
open neighbourhood of $\widehat\rho$, we obtain a branch
of constant Morse index $k$ since the index of a continuous
family of nondegenerate symmetric forms is locally constant.
Moreover,  $U'\ \subset \ U$, so the preserved negative subspaces
give $k \, \geq \,  r$.

Finally,  $V'$ is nonempty and open and by Theorem~\ref{thm:morse-energies} 
$$\rho^\dagger \in  V' \cap \mathcal R_n. $$
\end{proof}
An infinite Morse critical set is still possible under this argument.
If it occurs,  it is countable and every sequence of distinct critical
points has energy tending to infinity. This follows by exhausting
$\T$ with the sublevels $\bigl\{E_\rho\, \leq\, a\bigr\}$, 
$a\,  \in \, \mathbb{N}$.

\section{From instability to an open family}\label{sec:stable}

\subsection{Morse parity and the index bound}\label{subsec:morse-inputs}

Now, we prove persistence of critical points near an unstable pair. For this, we recall the following standard consequence of Morse theory, which we will use in the following proof~\cite[Theorem~3.5]{Milnor}. Let $f \colon  M\,  \longrightarrow \, \R$ be a smooth,  proper, 
bounded-below Morse function on a manifold diffeomorphic to $\R^d$.
If its critical set is finite,  then
\begin{equation}\label{eq:parity}
                  \sum_{x\,  \in \, \Crit(f)}{(-1)}^{ \ind_x(f)}\,  = \, 1.
\end{equation}
Therefore, if $f$ has a positive-index critical point,  it has at
least three distinct critical points, whether its critical set is
finite or infinite.
\begin{theorem}[Theorem~\ref{intro:persistence}]\label{thm:local-persistence}
Let $g\, \geq \, 2$ and $n\, \geq \, 2$,  let $\rho_*\,  \in \, \Hit$,  and suppose
$J_*\,  \in \, \Crit(E_{\rho_*})$. Assume that
$\Hess_{J_*}E_{\rho_*}$ has a negative direction; its kernel
may be nontrivial. For every open neighbourhood $V\, \subset\, \Hit$
of $\rho_*$ and every open neighbourhood $O\, \subset\, \T$ of $J_*$, 
there exist a nonempty connected open set $W\, \subset\, V$ and
real-analytic maps
$$
 J_0, \, J_1, \, J_2 \colon  W\,  \longrightarrow \, \T
$$
such that,  for every $\rho\,  \in \, W$,  their values are pairwise
distinct nondegenerate critical points of $E_\rho$. 
$J_0(\rho)$ is a strict local minimum.  Moreover,  $J_1(W)\, \subset\, O$, 
and $J_1(\rho)$ has constant Morse index $k$ on $W$,  where
$1\, \leq\, k\, \leq\, 3g-3$. If the Hessian is negative definite
on an $r$-dimensional real subspace,  the branch $J_1$ can be
chosen with $k\, \geq \, r$.
\end{theorem}
\begin{proof}
For the first statement, we choose any one-dimensional subspace where the Hessian is negative definite.  For the last assertion,  we use the given $r$-dimensional subspace. By Proposition~\ref{prop:local-regularisation}, we have a connected open
set $V'\, \subset\, V$,  and an analytic nondegenerate critical branch
$J_1 \colon  V'\,  \longrightarrow \, O$ of constant index $k\, \geq \, r$,  and
a representation $\rho^\dagger\,  \in \, V' \cap \mathcal R_n$.

$E_{\rho^\dagger}$ is Morse. It has the
positive index at $J_1(\rho^\dagger)$.  Since $\T\, \cong\, \R^{6g-6}$ and $E_{\rho^\dagger}$ is proper, we get at least three distinct critical points as discussed earlier. We choose a global minimum $J_0^\dagger$ and a third point $J_2^\dagger$, in addition to $J_1(\rho^\dagger)$. All three
are nondegenerate.

Choose disjoint coordinate neighbourhoods of these three points. Hessians in these three points are nondegenerate,  so Lemma~\ref{lem:regular}
and~\eqref{eq:hessian} show that $D_J\Phi$ is invertible at each selected
pair.  By the implicit function theorem, we get  one critical branch in each neighbourhood for all
representations in a sufficiently small common neighbourhood of
$\rho^\dagger$. Choose a
connected open $W\, \subset\, V'$ in that common neighbourhood. The newly
chosen branches are denoted $J_0$ and $J_2$,  and the previous $J_1$
is restricted to $W$. Their images remain in the disjoint neighbourhoods, 
so they are pairwise distinct. Their Hessians remain invertible.
The index-zero branch $J_0$ remains a strict local minimum,  and
$J_1(W)\, \subset\, O$ still has index $k\, \geq \, r$. The upper bound $3g-3$
follows from~\cite[Corollary~4.1]{SlegersPSH},  which states that the Morse index
of $J_1(\rho)$ is at most $3g-3$.
\end{proof}
\subsection{From area instability to energy instability}\label{subsec:area-to-energy}
 Let $u_0$ be a minimal immersion for a fixed representation, 
and suppose that a smooth equivariant variation $u_t$ satisfies
$$
 \left.\frac{d}{dt}\right|_0\Area(u_t)\,  = \, 0, 
 \qquad
 \left.\frac{d^2}{dt^2}\right|_0\Area(u_t)\, <\, 0.
$$

\begin{lemma}\label{lem:area}
Suppose a $\rho$-equivariant minimal immersion $u_0$ has a variation
with negative second variation of area. If $J_0$ is its induced marked
conformal structure,  then $\Hess_{J_0}E_\rho$ has a negative direction.
\end{lemma}
\begin{proof}
The equivariant area index of $u_0$ equals the Morse index
of $E_\rho$ at $J_0$; see \cite[Section~4.2]{SS}
and the reference to Ejiri's index comparison therein.
Thus a negative area variation implies that
$\Hess_{J_0}E_\rho$ has a negative direction.
\end{proof}
\begin{corollary}\label{cor:open-families}
Let $g\, \geq \, 3$ and $n\, \geq \, 4$. For every Hitchin component $\Hit$, 
there exist a nonempty connected open set $W\, \subset\, \Hit$ and
real-analytic
maps
$$J_0, \, J_1, \, J_2 \colon  W\,  \longrightarrow \, \T $$
whose values are pairwise distinct nondegenerate critical points
of $E_\rho$ for every $\rho\,  \in \, W$.  $J_0(\rho)$ is
a strict local minimum,  and $J_1(\rho)$ has a constant Morse index on $W$. In particular, every representation in $W$
has at least three nondegenerate marked critical points of its energy.
\end{corollary}
\begin{proof}
    For $g\, \geq \, 3$ and $n\, \geq \,  4$,  we take the unstable minimal map as in  ~\cite[Theorem A]{SS}.   This map is an immersion by~\cite[Remark 4.3]{SS},  so Lemma~\ref{lem:area} and Theorem~\ref{thm:local-persistence} give the three branches.
\end{proof}
\subsection{Bounded-energy continuation}\label{sec:compact}
Sagman--Smillie prove that,  for every $\rho_0\,  \in  \, \Hit$ and $B\,  \in  \,  \R$,  there are a neighbourhood $U$ of $\rho_0$ and a compact set $K_B\, \subset\, \T$ such that 
\begin{equation}\label{eq:uniformproper}
    \bigcup_{\rho\,  \in \,  U} \bigl\{J\,  \in \, \T\, \big|\, E_\rho(J)\, \leq\, B\bigr\}
    \, \ \subset \ \, K_B.
\end{equation}
For this we refer to \cite[Definition 4.11 and proof of Lemma 4.10]{SS}.

We claim the following:
\begin{proposition}\label{prop:bounded-continuation}
    Let $\rho_0\,  \in \, \mathcal R_n$,  and let $B\,  \in \, \R$ not be a critical value of $E_{\rho_0}$. Then there exists a connected open neighbourhood $U\, \ \subset \ \, \Hit$ of $\rho_0$,  an integer $N\, \geq \, 0$,  and real-analytic maps
    $$J^i \colon  U\,   \longrightarrow  \, \T, \qquad 1\, \leq\, i\, \leq\, N$$
    such that,  for every $\rho\,  \in  \, U$,  
$$\bigl\{J\,  \in \, \Crit(E_\rho)\, \big|\, 
 E_\rho(J)\, \leq\, B\bigr\}
 \,  = \, \bigl\{J^i(\rho)\, \big|\, 1\, \leq\, i\, \leq\, N\bigr\}.$$
The branches are pairwise distinct and nondegenerate,  and their Morse
indices are constant on $U$.
\end{proposition}

\begin{proof}
The set
$$\Crit(E_{\rho_0}) \cap \bigl\{E_{\rho_0}\, \leq\, B\bigr\}
$$
is finite because the sublevel is compact and $E_{\rho_0}$ is Morse.
Since $B$ is not a critical value, we enumerate this set as
$J^1_0, \, \ldots, \, J^N_0$,  with every critical energy strictly below
$B$. Now choose pairwise disjoint neighbourhoods so that  the analytic
implicit function theorem gives unique critical branches
$J^i(\rho)$ in those neighborhoods. We choose neighborhoods such that these
branches remain pairwise distinct and nondegenerate,  their energies
remain below $B$,  and their indices are constant.

Apply~\eqref{eq:uniformproper} at level $B$. If these branches did not
exhaust the critical points of energy at most $B$ on every sufficiently
small neighbourhood of $\rho_0$,  there would be sequences
$$
 \rho_a\,  \longrightarrow \, \rho_0, \qquad
 J_a\,  \in \, \Crit(E_{\rho_a}), \qquad E_{\rho_a}(J_a)\, \leq\, B, 
$$
with $J_a$ outside the chosen uniqueness neighbourhoods. Passing to a subsequence in the common compact set $K_B$, we have 
$J_a\,  \longrightarrow \, J_ \infty$. Since $E$ and $\Phi$ are jointly continuous, we have 
$$J_ \infty\,  \in \, \Crit(E_{\rho_0}), \qquad
 E_{\rho_0}(J_ \infty)\, \leq\, B.
$$
Thus $J_ \infty\,  = \, J^i_0$ for some $i$. Since $J^i_0$ lies in the
interior of its uniqueness neighbourhood,  $J_a$ must eventually lie in
that neighbourhood,  contradicting the choice of $J_a$.
\end{proof}


\bigskip 
\bibliographystyle{alpha}
\bibliography{references}
\end{document}